\documentclass[12pt]{amsart}

\usepackage{graphicx}

\usepackage{amsmath, mathrsfs, amssymb, eucal, amsthm, amscd, amsxtra, enumerate,faktor,enumitem}

\usepackage{geometry}

\usepackage[toc,page]{appendix}
\usepackage{array}
\usepackage{diagbox}
\usepackage{booktabs}
\usepackage{siunitx}
\usepackage{float}
\restylefloat{table}
\usepackage{caption}
\usepackage{placeins}
\usepackage{microtype}
\usepackage{hyperref}

\hypersetup{
colorlinks=true,
linkcolor=blue,
citecolor=red,
filecolor=magenta,      
urlcolor=cyan,
pdftitle={Mordell_Curves_Arithmetic_Progression},
pdfpagemode=FullScreen,
}

\usepackage{amsrefs}

\usepackage{cleveref}

\numberwithin{equation}{section}
\newtheorem{theorem}{\bf Theorem}[section]
\newtheorem{lemma}[theorem]{\bf Lemma}

\theoremstyle{remark}

\title{Powers as Fibonacci Sums}

\author[B. Earp-Lynch]{Benjamin Earp-Lynch}
\author[S. Earp-Lynch]{Simon Earp-Lynch}
\address{Department of Mathematics, Carleton University, Canada}
\email{benjaminearplynch@cmail.carleton.ca}
\email{simonearplynch@cmail.carleton.ca}
\author[O. Kihel]{Omar Kihel}
\address{Department of Mathematics and Statistics, Brock University, Canada}
\email{okihel@brocku.ca}
\author[P. Tiebekabe]{Pagdame Tiebekabe}
\address{Facult\'e des Sciences et Techniques, D\'epartement de Math\'ematiques, Universit\'e de Kara, Togo and Laboratoire d'Alg\`ebre, de Cryptologie, de G\'eom\'etrie Alg\'ebrique et Applications (LACGAA), Universit\'e Cheikh Anta Diop de Dakar (UCAD), S\'en\'egal}
\email{pagdame.tiebekabe@ucad.edu.sn}
\begin{document}
\begin{abstract}
We examine the equation $y^{a}=\sum\limits_{i=1}^{k}F_{n_{i}}$ for positive integers $y,a\geq 2$ and $k\geq3$.  This equation can be expressed as a problem in terms of the Zeckendorf representations of integers.  Using bounds on linear forms in logarithms and Baker-Davenport reduction methods, we are able to completely solve the equation for $ y\leq 25000$ when $k=3$, for $y\leq 1000$ when $k= 4$, for $y\leq 40$ when $k= 5$, and for $y\leq 3$ when $k=6$.
\end{abstract}
\maketitle
\section{Introduction}
Let $y$ and $a\geq 2$ be positive integers, and let $F_{r}$ denote the $r$-th Fibonacci number.   The equation 
\begin{equation}\label{Main}
y^{a}=\sum\limits_{i=1}^{k}F_{n_{i}},
\end{equation}
has been solved when $k=1$ for any $y$ by Bugeaud, Mignotte and Siksek in \cite{MR2215137}, and for $(y,k)\in \{(2,2),(2,3),(2,4),(2,5)\}$ 
respectively by Bravo and Luca in \cite{MR3510734}, by Bravo and Bravo in \cite{MR3379028}, and by Tiebekabe and Diouf in \cite{MR4370726} and \cite{MR4539734}.  For $k=2$ and $y$ an arbitrary positive integer greater than or equal to $2$,  Kebli, Kihel, Larone and Luca in \cite{MR4177538} showed that the equation has finitely many solutions, and solved it for $y\leq 1000$.  In \cite{MR4316638}, again for $k=2$, Kihel and Larone found a bound for $a$ in terms of $y$. In \cite{MR4623323}, Ziegler found that in the case when $k=2$ and $a\geq 2$, there is at most one solution except when $y=2,3,4,6,10$, and in \cite{MR4246600}, Vukusic and Ziegler found an upper bound for $y^{a}$ in the case when $k=2$.  The upper bound for $y^a$ in the case when $k=2$ in \cite{MR4246600} depended on the Hamming weight of the Zeckendorf representation of $y$.  

Let $N\in \mathbb{Z}^{+}$. Then there is a unique sum $$N=\epsilon_{l}F_{l}+\epsilon_{l-1}F_{l-1}+\dots +\epsilon_{2}F_2,$$ with $\epsilon_{l}=1$, $\epsilon_{j}\in \{0,1\}$ and $\epsilon_{j}\epsilon_{j+1}=0$ for all $j \in \{2,\dots l-1\}$.  This is called the Zeckendorf representation of $N$, as first described by Zeckendorf in \cite{MR0308032}. The number of positive integers $j$ for which $\epsilon_{j}$ is equal to $1$ is referred to as the number of digits in the Zeckendorf representation of $N$. 

For a given $k$, the problem of representing powers of $y$ as a sum of $k$ Fibonacci numbers may be reduced to that of finding those powers of $y$ with no more than $k$ digits in their Zeckendorf representation. 
In this paper, by considering the sum on the right hand side of \Cref{Main} as a Zeckendorf representation, we extend results from \cites{MR3510734, MR3379028, MR4370726, MR4539734, MR4177538}, solving \eqref{Main} completely for $y\leq25000$ when $k=3$, for $y\leq 1000$ when $k\leq 4$, for $y\leq 40$ when $k\leq 5$, and for $y\leq3$ when $k=6$. 
The following theorem contains a summary of our results, referring to the solutions which can be found in the tables compiled in \Cref{SolutionsSection}.
\begin{theorem}\label{MainResult}
Let $y$ and $a$ be integers greater than or equal to $2$. 
\begin{itemize}
    \item For $k=3$ and $y\leq 25000$, all solutions to \eqref{Main} are contained in Table \ref{tab:k3_solutions}.
    \item For $k=4$ and $y\leq 1000$, all solutions to \eqref{Main} are contained in Table \ref{tab:k4_solutions}.
    \item For $k=5$ and $y\leq 40$, all solutions to \eqref{Main} are contained in Table \ref{tab:k5_solutions}.
    \item For $k=6$ and $y\leq 3$, all solutions to \eqref{Main} are contained in Table \ref{tab:k6_solutions}.
\end{itemize}
\end{theorem}

We obtain a bound on the size of the largest digit in the Zeckendorf representation of a power $y^a$ with $k$ Zeckendorf digits.  Our setup, in order to apply a theorem of Matveev to bound a linear form in three logarithms, is similar to that used by Bugeaud in the proof of Lemma 3.3 in \cite{MR4219035}.  In that work it was shown that, for a given $k$, there are finitely many integers with $k$ or fewer digits in their Zeckendorf representation and whose prime divisors belong to a given finite set.  In contrast to the work in \cite{MR4219035}, we compute the constants explicitly.  We apply the bound obtained to the case of finding all powers of $y$ with $k$ or fewer Zeckendorf digits. We then invoke multiple applications of the Baker-Davenport reduction method, wherein bounds on parameters are applied to previous steps in order to tighten the bounds as the procedure goes along. This significantly reduces the computation required by the time one arrives at the later reductions. Finally, we are able to find all solutions to \eqref{Main} for the values of $y$ and $k$ mentioned above. In the case $k=3$, these include those belonging to the infinite family $$L_{2s}^{2}=F_{4s+1}+F_{4s-1}+F_{3},$$ where $L_{2s}$ denotes the $2s$th Lucas number. The existence of this family can be seen from the identities $L_{s}^{2}=L_{2s}+2(-1)^s$ and $L_{s}=F_{s+1}+F_{s-1}$. All computations were run using SageMath \cite{sagemath}, and the code can be found at the url listed in the references \cite{CodeLink}.

\section{Preliminaries}\label{Prelims}

Some previous work on \eqref{Main} has allowed for the sum to contain repeated or consecutive Fibonacci numbers.  However, it is sufficient to consider only Zeckendorf representations, since the Zeckendorf representation of an integer is unique, and any other representation of an integer as a sum of Fibonacci numbers can be turned into its Zeckendorf representation through repeated use of the identities $F_{m-1}+F_{m-2}=F_{m}$ and $2F_{m_{i}}=F_{m_{i}+1}+F_{m_{i}-2}$, so that the resulting Zeckendorf sum contains at most the number of terms in the original sum. 

To obtain the initial bounds, we closely follow the methods used in the proofs of Lemmas 3.2 and 3.3 in \cite{MR4219035}. Suppose that $y$ and $a$ are positive integers, and $y^{a}$ has Zeckendorf representation $$y^{a}=\sum\limits_{i=1}^{k}F_{n_{i}},$$ with $n_{1}\geq 2$ and $n_{i+1}\geq n_{i}+2$. Since $n_{k}-n_{i}\geq2(k-i)$ and $n_{k}+n_{i}\geq2(k+i)$, it follows from the triangle inequality that $$y^{a}\leq\frac{\alpha^{n_{k}}}{\sqrt{5}}\left(\sum\limits_{l=0}^{\infty}\frac{1}{\alpha^{2l}}\right)=\frac{\alpha^{n_{k}+1}}{\sqrt{5}}<\alpha^{n_{k}},$$ and so $a\log{y}<n_{k}\log{\alpha}$.

Let $\gamma$ be an algebraic number with minimal polynomial $$f(x)=a_{0}\prod\limits_{i=1}^{d}\left(x-\gamma^{(i)}\right),$$ then the absolute logarithmic height of $\gamma$ is defined $$h(\gamma)=\frac{1}{d}\left(\log{a_{0}}+\sum\limits_{i=1}^{d}\log\left(\max\left\{\left|\gamma^{(i)}\right|,1\right\}\right)\right).$$ The following three properties of the absolute logarithmic height are well-known.
\begin{align*}
&h(\gamma+\delta)\leq h(\gamma)+h(\delta)+\log2,\\
&h(\gamma\delta^{\pm1})\leq h(\gamma)+h(\delta),\\
&h(\gamma^{k})=|k|h(\gamma).
\end{align*} 

Using these properties, we find an upper bound on the absolute logarithmic height of the quantity $\alpha^{n_{k}-n_{j}}+\dots+\alpha^{n_{j+1}-n_{j}}+1$. We have
\begin{align*}
h\left(\alpha^{n_{k}-n_{j}}+\dots+\alpha^{n_{j+1}-n_{j}}+1\right)\leq& (k-j)\log2+(k-j)\left(n_{k}-n_{j}\right)\frac{\log\alpha}{2}\\
\leq& \left(\frac{k}{2}\log{2\alpha}\right)\left(n_{k}-n_{j}\right).
\end{align*}

For $2\leq j\leq k$, set
\begin{align*}
\Lambda_{j}:=&\left|y^{a}\alpha^{-n_{j}}\sqrt{5}\frac{1}{\alpha^{n_{k}-n_{j}}+\alpha^{n_{k-1}-n_{j}}+\dots+1}-1\right|\\
=&\frac{\sum\limits_{h=1}^{j-1}\sqrt{5}F_{n_{h}}-\sum\limits_{h=j}^{k}\beta^{n_{h}}}{\sum\limits_{h=j}^{k}\alpha^{n_{h}}}
\leq\alpha^{n_{j-1}-n_{k}}\left(\sum\limits_{h=1}^{j-1}\frac{1}{\alpha^{n_{j-1}-n_{h}}}+\sum\limits_{h=1}^{k}\frac{1}{\alpha^{n_{j-1}+n_{h}}}\right)\\
\leq&\alpha^{n_{j-1}-n_{k}}\left(\sum\limits_{l=0}^{\infty}\frac{1}{\alpha^{2l}}\right)=\alpha^{1+n_{j-1}-n_{k}}.
\end{align*}

The last inequality follows from the fact that this is a Zeckendorf representation, and so $n_{j-1}-n_{h}\geq 2(j-1-h)$ and $n_{j-1}+n_{h}\geq 2(j-1+h)$.  In particular, this gives us the inequality
\begin{equation}\label{ineq1}
\Lambda_{j}\leq \frac{\alpha}{\alpha^{n_{k}-n_{j-1}}}.
\end{equation}

Before we proceed, it will be necessary to establish that $\Lambda_{j}$ never vanishes. To see this, observe that $\Lambda_{j}=0$ would mean $$y^{a}\sqrt{5}=\sum\limits_{i=j}^{k}\alpha^{n_{i}}.$$ Using the fact that $\alpha^{\ell}=\frac{L_{\ell}+F_{\ell}\sqrt{5}}{2}$ for any positive integer $\ell$, we would then have $$\left(2y^{a}-\sum\limits_{i=j}^{k}F_{n_{i}}\right)\sqrt{5}=\sum\limits_{i=j}^{k}L_{n_i}.$$ The left side of this equation is either $0$ or irrational, whereas the right side is always a positive integer, and so the equality can never hold.

It now follows from \eqref{ineq1} that $$\log\Lambda_{j}\leq-(n_{k}-n_{j-1}-1)\log\alpha.$$ 

\section{Initial Bounds}
We will use the following well known theorem of Matveev \cite{MR1817252}, from which we can acquire a lower bound for $\log \Lambda_{j}$.

\begin{theorem}\label{Matveev}
Let $n\geq 2$ be an integer. Let $\gamma_{1}, \dots, \gamma_{n}$ be non-zero algebraic real numbers. Let $D$ be the degree over $\mathbb{Q}$ of a number field containing $\gamma_{1}, \dots, \gamma_{n}$. Let $A_{1}, \dots, A_{n}$ be real numbers with $$\log A_{j} \geq \max \left\{h(\gamma_{j}), \frac{|\log \gamma_{j}|}{D},\frac{0.16}{D}\right\}, \quad 1\leq j \leq n.$$ Let $b_{1}, \dots, b_{n}$ be integers and set $$B'=\max\left\{1, \max\left\{|b_{j}|\frac{\log A_{j}}{\log A_{n}}: 1\leq j \leq n \right\} \right\}.$$
Then, provided that $\gamma_{1}^{b_{1}}\cdots\gamma_{n}^{b_{n}}\neq1$, we have $$\log\left|\gamma_{1}^{b_{1}}\cdots \gamma_{n}^{b_{n}}-1\right|>-2\times 30^{n+3}n^{4.5}D^{n+2}\log(eD) \log A_{1}\cdots \log A_{n} \log(eB').$$
\end{theorem}

We invoke \Cref{Matveev} with $n=3,D=2$ and
\begin{alignat*}{3}
\gamma_{1}&=y, &\quad \gamma_{2} &= \alpha, &\quad \gamma_{3} &= \frac{\alpha^{n_{k}-n_{j}}+\dots+\alpha^{n_{j+1}-n_{j}}+1}{\sqrt{5}},\\
b_{1} &= a, &\quad b_{2} &= -n_{j}, &\quad b_{3}&=-1.
\end{alignat*}
We may take $$A_{1}=y,\quad A_{2}=\alpha,\quad A_{3}=\exp\left( \frac{k}{2}\left(n_{k}-n_{j}\right)\log2\alpha+\log{5}\right)$$
and $$B'=\frac{n_{k}}{(n_{k}-n_{j})+\frac{\log{5}}{\frac{k}{2}\log{2\alpha}}},$$ which follows since $a\log{y}<n_{k}\log{\alpha}$.

\noindent Since we have shown that $\Lambda_{j}$ is never $0$, we have that $\gamma_{1}\gamma_{2}\gamma_{3}\neq1$ and \Cref{Matveev} yields the lower bound $$\log\Lambda_{j}>-C_{k}\left( \left(n_{k}-n_{j}\right)+m_{k}\right) \log\left(e\frac{n_{k}}{(n_{k}-n_{j})+m_{k}}\right),$$ where $C_{k}=2\times 30^{6}3^{4.5}2^{5}\log(2e) \log{y}\log{\alpha}\left(\frac{k}{2}\log{2\alpha}\right)$ and $m_{k}=\frac{\log{5}}{\frac{k}{2}\log{2\alpha}}$. 

\noindent From this, we obtain 
\begin{align*}
(n_{k}-n_{j-1}-1)\log\alpha<&C_{k}\left( \left(n_{k}-n_{j}\right)+m_{k}\right) \log\left(\frac{en_{k}}{(n_{k}-n_{j})+m_{k}}\right)\\
<&C_{k}\left( \left(n_{k}-n_{j}\right)+m_{k}\right) \log\left(\frac{en_{k}}{2}\right).
\end{align*}

\noindent Hence, 
\begin{equation}\label{ineq2}
n_{k}-n_{j-1}< \frac{C_{k}}{\log\alpha}(n_{k}-n_{j}+m_k)\log{\left(\frac{en_{k}}{2}\right)}+1
\end{equation}
for all $2\leq j\leq k$. 

Setting $$f(x):=\frac{C_{k}}{\log\alpha}(x+m_{k})\log{\left(\frac{en_{k}}{2}\right)}+1,$$ we see that 
\begin{equation}\label{k-1Bd}
n_{k}-n_{1}< f^{k-1}(0),
\end{equation}
where $f^{k-1}$ denotes the $k-1$st iterate of the function $f$.
Observe that
\begin{align*}
0<\left| \sqrt{5}y^{a}-\alpha^{n_{1}}-\dots-\alpha^{n_{k}}\right|\leq \left|\beta^{n_{1}}\right|+\dots+\left|\beta^{n_{k}}\right|<\alpha,
\end{align*}
where the first inequality can be shown in a similar manner to the demonstration at the end of \Cref{Prelims} that the $\Lambda_{j}$ are nonzero. and so
\begin{align*}
\Lambda=\left|\sqrt{5}y^{a}\alpha^{-n_{1}}\frac{1}{1+\alpha^{n_{2}-n_{1}}+\dots+\alpha^{n_{k}-n_{1}}}-1\right|<\alpha^{-n_{k}+1}.
\end{align*}

We apply \Cref{Matveev} again, with constants $m=3,D=2$ and

\begin{alignat*}{3}
\gamma_{1} &= y, &\quad \gamma_{2} &= \alpha, &\quad \gamma_{3} &= \frac{\sqrt{5}}{1+\alpha^{n_{2}-n_{1}}+\dots+\alpha^{n_{k}-n_{1}}},\\
b_{1} &= a, &\quad b_{2} &= -n_{1}, &\quad b_{3} &= 1,
\end{alignat*}
for which we may take $A_{1}=y$, $A_{2}=\alpha$ and $A_{3}\geq \exp\left((n_{k}-n_{1})\left(k\log{\alpha}+k+1\right)\right)$, with $$B'\geq\frac{\log\alpha}{k\log\alpha+k+1}n_{k}.$$ This gives the bound $$\log\Lambda>-C_{\Lambda}(n_{k}-n_{1})\log\left(e\frac{\log\alpha}{k\log\alpha+k+1}n_{k}\right),$$ where $$C_{\Lambda}=2\times 30^{6}3^{4.5}2^{5}\left(k\log{\alpha}+k+1\right)\log(2e) \log y\log \alpha,$$

\noindent and so from above we get 
\begin{equation}\label{ineq2.5}
(n_{k}-1)\log{\alpha}<C_{\Lambda}(n_{k}-n_{1})\log\left(e\frac{\log\alpha}{k\log\alpha+k+1}n_{k}\right).
\end{equation} 
Using the bound \eqref{k-1Bd} obtained for $n_{k}-n_{1}$, we see $$n_{k}<\frac{C_{\Lambda}}{\log\alpha}f^{k-1}(0)\log\left(e\frac{\log\alpha}{k\log\alpha+k+1}n_{k}\right)+1.$$ When $k$ and $y$ are fixed integers, the right side of the above expression is a polynomial of degree $k$ in $\log{n_{k}}$, and so it yields a bound on $n_{k}$. 

\section{Reducing the Bounds}

We will now outline how the Baker-Davenport reduction method may be applied in a way which will reduce the computation required for a complete solution significantly. First, we may obtain an upper bound on $n_{k}-n_{k-1}$, which leads to an improved bound on $n_{k}$.

Recall from \eqref{ineq1} that $\Lambda_{k}<\frac{\alpha}{\alpha^{n_{k}-n_{k-1}}},$ and defining $z_{1}:=\sqrt{5}p^{a}\alpha^{-n_{k}}$, we obtain
\begin{align*}
1<1+\alpha^{-n_{k}+1}&=1+\alpha^{-n_{k}}\left(\sqrt{5}+\left(\frac{1-\sqrt{5}}{2}\right)\right)\leq 1+\alpha^{-n_{k}}\left(\sqrt{5}-\beta^{n_{k}}\right)\\
&\leq 1-\alpha^{-n_{k}}\beta^{n_{k}}+\sqrt{5}\alpha^{-n_{k}}\left(F_{n_{k-1}}+\dots+F_{n_{1}}\right) =z_{1}. 
\end{align*}
It follows that $$0<\log z_{1}\leq e^{\log z_{1}}-1=\Lambda_{k}<\frac{\alpha}{\alpha^{n_{k}-n_{k-1}}}.$$ After division by $\log\alpha$, this gives $$0<\left|a\frac{\log{y}}{\log{\alpha}}-n_{k}+\frac{\log{\sqrt{5}}}{\log{\alpha}}\right|<\frac{\alpha}{\log\alpha}\alpha^{-(n_{k}-n_{k-1})}.$$ Now, we apply the following version of the Baker-Davenport reduction method, due to Bravo, G\'omez and Luca \cite{MR3527869}.

\begin{lemma}[Bravo, G\'omez and Luca]\label{Reduction}
Let $N$ be a positive integer. Let $p/q$ be a convergent of the continued fraction expansion of the irrational $\kappa$ such that $q>6N$, and let $A,B,\mu$ be real numbers with $A>0$ and $B>1$. Furthermore, let $\varepsilon=||\mu q||-N\cdot||\kappa q||$, where $||\cdot||$ denotes the distance from the nearest integer. If $\varepsilon>0$, then there is no solution of the inequality $$0<\left|u\kappa-v+\mu\right|<AB^{-w}$$ in positive integers $u,v$ and $w$ with $$u\leq N\quad\text{and}\quad w\geq\frac{\log{(Aq/\varepsilon)}}{\log{B}}.$$
\end{lemma}

In the notation of the lemma, we will have 
\begin{alignat*}{3}
u &= a, &\quad \kappa &= \frac{\log{y}}{\log\alpha}, &\quad v &= n_{k}, \\
\mu &= \frac{\log{\sqrt{5}}}{\log\alpha}, &\quad A &= \frac{\alpha}{\log\alpha}, &\quad B &= \alpha, \\
w &= n_{k}-n_{k-1}. &&&&
\end{alignat*}
and the bound $N$ may be obtained as outlined above.

As a result of this reduction, we may obtain a bound $N_{k-1}$ on $n_{k}-n_{k-1}$. Examining the work from earlier, in particular inequality \eqref{ineq2}, we see that this allows us to obtain a better bound, $N^{(1)}$, for $n_{k}$ using \begin{equation}\label{reductionineq}
n_{k}<\frac{C_{\Lambda}}{\log\alpha}f^{k-1}(N_{k-1})\log\left(e\frac{\log\alpha}{k\log\alpha+k+1}n_{k}\right)+1.
\end{equation}

Using this bound, we can go through the reduction procedure again, this time setting $N=N^{(1)}$, to obtain the bound $n_{k}-n_{k-1}\leq N_{k-1}^{(1)}$, which may itself allow us to reduce the bound on $n_{k}$ slightly further, but we found that little improvement was made when carrying out the reduction a third time at this point. This is due to the fact that the next bound on $n_{k}$ was usually not enough smaller to allow the use of a smaller value of $q$.%

In general, using inequality \eqref{ineq1}, and given a bound $N$ on $n_{k}$ and bounds $N_{k-1},\dots,N_{j}$ on $n_{k}-n_{k-1},\dots,n_{k}-n_{j}$, we apply \Cref{Reduction} to the inequality 
\begin{equation}\label{ineq3}
0<\left|a\frac{\log y}{\log\alpha}-n_{k}+\frac{\log\Upsilon_{j}(t_{k-1},\dots,t_{j})}{\log\alpha}\right|<\frac{\alpha}{\log\alpha}\alpha^{-(n_k-n_{j-1})},
\end{equation}
where $\Upsilon_{j}(t_{k-1},\dots,t_{j})=\frac{\sqrt{5}}{1+\alpha^{t_{k-1}}+\dots+\alpha^{t_{j}}}$ for every possible combination of integers $t_{k-1},\dots,t_{j}$ with $t_{i}<N_{i}$ and $t_{l}\leq t_{m}-2$ for $l>m$. In the notation of \Cref{Reduction}, we have 
\begin{alignat*}{3}
u &= a, &\quad \kappa &= \frac{\log{y}}{\log\alpha}, &\quad v &= n_{k},\\
\mu &= \frac{\log{\Upsilon_{j}(t_{k-1},\dots,t_{j})}}{\log\alpha}, &\quad A &= \frac{\alpha}{\log\alpha}, &\quad B&=\alpha,\\ w &= n_{k}-n_{j-1}. &&&&
\end{alignat*} 
This gives a bound on $n_{k}-n_{j-1}$, and when $j=1$,  \eqref{ineq2.5} gives a much lower bound on $n_{k}$ itself. 

Problems may arise, however, when $\mu=s_{1}\frac{\log{\Upsilon_{j}(t_{k-1},\dots,t_{j})}}{\log\alpha}=s_{2}\frac{\log y}{\log\alpha}+s_{3}$ for integers $s_{1},s_{2}$ and $s_{3}$, since this may make it impossible to find a suitable convergent with denominator $q$ that makes $\varepsilon=||\mu q||-N\cdot||\kappa q||$ positive as required by \Cref{Reduction}. 
Note that we will abbreviate $\Upsilon_{j}(t_{k-1},\dots,t_{j})$ by simply $\Upsilon_{j}$ going forward, and include the $t_{i}$ when they are are once more relevant. 
 When $$s_{1}\log\Upsilon_{j}=s_{2}\log y+s_{3}\log\alpha,$$ we have $\Upsilon^{s_{1}}=y^{s_{2}}\alpha^{s_{3}}$. Taking the absolute values of the norms of both sides in $\mathbb{Q}(\sqrt{5})$, we have $|N(\Upsilon_{j})|^{s_{1}}=y^{2s_{2}}$, and so we obtain $\frac{2s_{2}}{s_{1}}=\frac{\log{|N(\Upsilon_{j})|}}{\log{y}}$. The condition $\frac{\log{|N(\Upsilon_{j})|}}{\log{y}}\in\mathbb{Q}$ may be checked in SageMath \cite{sagemath} for every possible value of $\Upsilon_{j}$, after which a suitable integer $s_{3}$ may be found if it exists. Subsequently, we may write inequality \eqref{ineq3} as $$\left|(as_{1}+s_{2})\frac{\log{y}}{\log{\alpha}}-(n_{k}s_{1}-s_{3})\right|<\frac{w\alpha}{\log{\alpha}}\alpha^{-(n_{k}-n_{j})}.$$ Let $[a_{0};a_{1},a_{2},\dots]$ be the continued fraction expansion of $\frac{\log{y}}{\log{\alpha}}$, with $p_{i}/q_{i}$ denoting its $i$th convergent. Since $a\leq n_{k}<N$, it follows that $as_{1}+s_{2}<Ns_{1}+s_{2}$. We compute this bound, find the convergent $p_{l}/q_{l}$ such that $q_{l}<Nw+r<q_{l+1}$ and note that the properties of continued fractions also give us the bound $$\left|(as_{1}+s_{2})\frac{\log{y}}{\log{\alpha}}-(n_{k}s_{1}+s_{3})\right|>\frac{1}{(a_{M}+2)(as_{1}+s_{2})},$$ where $a_{M}:=\max\{a_{i}:1\leq i\leq l\}$. Combining these two bounds gives us $$n_{k}-n_{j}<\frac{\log{\frac{(a_{M}+2)(Ns_{1}+s_{2})w\alpha}{\log{\alpha}}}}{\log{\alpha}}.$$ We compare this with the current bound on $N_{j}$ for the cases where there is no linear dependence and take the maximum.

We wrote a program \cite{CodeLink} in SageMath \cite{sagemath} to carry out the procedure of reducing the bounds and determining all possible solutions 
for given inputs $k$ and $y$ and subsequently to check for solutions to \eqref{Main} with $n_{k}$ up to the modest final bound. 
Our program can accommodate higher bounds on $k$ and $y$ than those listed here, which were chosen to be demonstrative while limiting the runtime. 

We exhibit the constants used to find all solutions in the case when $k=6$ and $y=2$.  First, using \Cref{Matveev}, we find an initial upper bound on $n_6$ of $5.5\times 10^{94}$.  Next, using the reduction method of \Cref{Reduction}, we are able to obtain an upper bound on $n_6-n_5$ of $464$. We substitute this into inequality \ref{reductionineq}, which in turn leads to an improved bound on $n_6$, which we can use to again improve the bound on $n_6-n_5$ to $404$, and this again improves the upper bound on $n_6$ to $3.9 \times 10^{81}$.  Subsequent iterations of this procedure do not prove useful.  However we then apply \Cref{Reduction} to obtain a bound on $n_6-n_4$ of $413$, and using this bound we again are able to further improve our bound on $n_k$, and subsequently we repeat our reductions to improve our bounds on $n_6-n_5$ and $n_6-n_4$.  At this step in the procedure, we have obtained an upper bound of $2.6 \times 10^{65}$ on $n_6$ and upper bounds of $324$ and $336$ on $n_6-n_5$ and $n_6-n_4$ respectively.  Proceeding in this manner, we are able to improve the bound on $n_6$ to $2.3 \times 10^{49}$, then to $2.5 \times 10^{33}$, then $4.1\times 10^{17}$, and finally to $n_6\leq 51$, with final bounds on $n_6-n_j$ of $47, 41, 33, 26, 23$ for $j=1, \dots 5$ respectively.  At this point, it is computationally manageable to find all possible solutions to \eqref{Main} with $k=6$ and $y=2$.

Below we include tables containing all solutions we computed to equation \eqref{Main}.

\section{Tables of Solutions}\label{SolutionsSection}
The tables in this section contain the solutions to Equation \eqref{Main}, up to the bounds listed in \Cref{MainResult}.

\begin{table}[htbp]
    \centering
    \small 
    \setlength{\tabcolsep}{4pt} 
    \caption{Solutions in the case $k = 3$.}
    \label{tab:k3_solutions}

    \begin{tabular}{@{}ccccc @{\hspace{2em}}|@{\hspace{2em}} ccccc@{}}
        \toprule
        $y$ & $a$ & $n_1$ & $n_2$ & $n_3$ & $y$ & $a$ & $n_1$ & $n_2$ & $n_3$ \\
        \midrule
        $L_{2s}$ & $2$ & $3$ & $4s-1$ & $4s+1$ & $17$ & $2$ & $2$ & $10$ & $13$ \\
        $2$ & $5$ & $4$ & $6$ & $8$          & $19$ & $3$ & $5$ & $11$ & $20$ \\
        $2$ & $6$ & $2$ & $6$ & $10$         & $20$ & $2$ & $3$ & $8$ & $14$ \\
        $2$ & $7$ & $5$ & $9$ & $11$         & $23$ & $2$ & $6$ & $12$ & $14$ \\
        $2$ & $8$ & $3$ & $8$ & $13$         & $24$ & $2$ & $10$ & $12$ & $14$ \\
        $2$ & $10$ & $4$ & $9$ & $16$        & $33$ & $2$ & $7$ & $11$ & $16$ \\
        $3$ & $3$ & $2$ & $5$ & $8$          & $35$ & $2$ & $5$ & $13$ & $16$ \\
        $3$ & $4$ & $5$ & $8$ & $10$         & $37$ & $2$ & $5$ & $14$ & $16$ \\
        $3$ & $5$ & $3$ & $6$ & $13$         & $57$ & $2$ & $10$ & $15$ & $18$ \\
        $5$ & $2$ & $2$ & $4$ & $8$          & $112$ & $2$ & $2$ & $17$ & $21$ \\
        $5$ & $3$ & $3$ & $9$ & $11$         & $139$ & $2$ & $7$ & $17$ & $22$ \\
        $5$ & $4$ & $3$ & $7$ & $15$         & $1161$ & $2$ & $10$ & $17$ & $31$ \\
        $6$ & $6$ & $10$ & $13$ & $24$       & $1476$ & $2$ & $9$ & $13$ & $32$ \\
        $7$ & $3$ & $8$ & $11$ & $13$        & $4999$ & $2$ & $12$ & $30$ & $37$ \\
        $10$ & $2$ & $4$ & $6$ & $11$        & & & & & \\
        \bottomrule
    \end{tabular}
\end{table}
%%%%%%%%%%%%

\newpage
\begin{table}[htbp]
    \centering
    \small 
    \setlength{\tabcolsep}{4pt} 
    \caption{Solutions in the case $k = 4$.}
    \label{tab:k4_solutions}

    \begin{tabular}{@{}cccccc @{\hspace{1.5em}}|@{\hspace{1.5em}} cccccc@{}}
        \toprule
        $y$ & $a$ & $n_1$ & $n_2$ & $n_3$ & $n_4$ & $y$ & $a$ & $n_1$ & $n_2$ & $n_3$ & $n_4$ \\
        \midrule
        $6$ & $4$ & $8$ & $10$ & $13$ & $16$    & $51$ & $2$ & $2$ & $4$ & $7$ & $18$ \\
        $6$ & $5$ & $4$ & $8$ & $16$ & $20$     & $60$ & $2$ & $6$ & $8$ & $16$ & $18$ \\
        $11$ & $2$ & $4$ & $6$ & $8$ & $11$     & $63$ & $2$ & $8$ & $14$ & $16$ & $18$ \\
        $11$ & $5$ & $10$ & $21$ & $23$ & $26$  & $65$ & $2$ & $3$ & $6$ & $9$ & $19$ \\
        $12$ & $3$ & $6$ & $9$ & $11$ & $17$    & $72$ & $2$ & $4$ & $7$ & $16$ & $19$ \\
        $13$ & $2$ & $2$ & $4$ & $8$ & $12$     & $83$ & $2$ & $2$ & $9$ & $11$ & $20$ \\
        $14$ & $2$ & $5$ & $7$ & $9$ & $12$     & $84$ & $2$ & $4$ & $10$ & $13$ & $20$ \\
        $14$ & $3$ & $4$ & $7$ & $12$ & $18$    & $92$ & $2$ & $7$ & $11$ & $17$ & $20$ \\
        $15$ & $2$ & $5$ & $8$ & $10$ & $12$    & $97$ & $2$ & $5$ & $10$ & $18$ & $20$ \\
        $15$ & $4$ & $8$ & $10$ & $19$ & $24$   & $105$ & $2$ & $4$ & $8$ & $10$ & $21$ \\
        $19$ & $2$ & $5$ & $9$ & $11$ & $13$    & $106$ & $2$ & $3$ & $10$ & $13$ & $21$ \\
        $21$ & $2$ & $2$ & $6$ & $10$ & $14$    & $119$ & $2$ & $8$ & $15$ & $18$ & $21$ \\
        $22$ & $2$ & $5$ & $7$ & $11$ & $14$    & $146$ & $2$ & $9$ & $16$ & $18$ & $22$ \\
        $23$ & $3$ & $2$ & $13$ & $16$ & $21$   & $170$ & $2$ & $3$ & $6$ & $13$ & $23$ \\
        $26$ & $2$ & $4$ & $6$ & $10$ & $15$    & $174$ & $2$ & $2$ & $8$ & $17$ & $23$ \\
        $28$ & $3$ & $5$ & $10$ & $19$ & $22$   & $358$ & $2$ & $2$ & $5$ & $20$ & $26$ \\
        $30$ & $2$ & $3$ & $10$ & $13$ & $15$   & $364$ & $2$ & $7$ & $12$ & $21$ & $26$ \\
        $31$ & $3$ & $4$ & $12$ & $16$ & $23$   & $382$ & $2$ & $10$ & $20$ & $22$ & $26$ \\
        $39$ & $2$ & $7$ & $12$ & $14$ & $16$   & $396$ & $2$ & $2$ & $20$ & $23$ & $26$ \\
        $41$ & $2$ & $6$ & $8$ & $10$ & $17$    & $445$ & $2$ & $3$ & $6$ & $17$ & $27$ \\
        $42$ & $2$ & $3$ & $8$ & $12$ & $17$    & $460$ & $2$ & $10$ & $19$ & $21$ & $27$ \\
        $48$ & $2$ & $6$ & $11$ & $15$ & $17$   & & & & & & \\
        \bottomrule
    \end{tabular}
\end{table}

\begin{table}[htbp]
    \centering
    \small
    \setlength{\tabcolsep}{5pt}
    \caption{Solutions in the case $k = 5$.}
    \label{tab:k5_solutions}
    \begin{tabular}{@{}ccccccc@{}}
        \toprule
        $y$ & $a$ & $n_1$ & $n_2$ & $n_3$ & $n_4$ & $n_5$ \\
        \midrule
        $2$ & $13$ & $6$ & $10$ & $14$ & $16$ & $20$ \\
        $3$ & $6$ & $2$ & $6$ & $8$ & $11$ & $15$ \\
        $6$ & $3$ & $2$ & $4$ & $7$ & $10$ & $12$ \\
        $10$ & $6$ & $10$ & $12$ & $24$ & $26$ & $30$ \\
        $11$ & $3$ & $2$ & $8$ & $11$ & $13$ & $16$ \\
        $15$ & $3$ & $4$ & $9$ & $12$ & $15$ & $18$ \\
        $20$ & $3$ & $3$ & $7$ & $13$ & $16$ & $20$ \\
        $28$ & $2$ & $2$ & $6$ & $8$ & $12$ & $15$ \\
        $31$ & $2$ & $6$ & $8$ & $11$ & $13$ & $15$ \\
        $33$ & $5$ & $7$ & $13$ & $15$ & $24$ & $38$ \\
        $34$ & $2$ & $2$ & $4$ & $8$ & $12$ & $16$ \\
        $37$ & $3$ & $3$ & $7$ & $11$ & $19$ & $24$ \\
        \bottomrule
    \end{tabular}
\end{table}

\begin{table}[htbp]
    \centering
    \small
    \setlength{\tabcolsep}{5pt}
    \caption{Solutions in the case $k = 6$.}
    \label{tab:k6_solutions}
    \begin{tabular}{@{}cccccccc@{}}
        \toprule
        $y$ & $a$ & $n_1$ & $n_2$ & $n_3$ & $n_4$ & $n_5$ & $n_6$ \\
        \midrule
        $2$ & $9$ & $2$ & $4$ & $6$ & $9$ & $11$ & $14$ \\
        $2$ & $11$ & $2$ & $5$ & $7$ & $10$ & $14$ & $17$ \\
        $2$ & $12$ & $2$ & $4$ & $12$ & $14$ & $16$ & $18$ \\
        $2$ & $14$ & $4$ & $9$ & $13$ & $16$ & $19$ & $21$ \\
        $2$ & $17$ & $6$ & $11$ & $13$ & $18$ & $20$ & $26$ \\
        $2$ & $28$ & $2$ & $8$ & $12$ & $20$ & $29$ & $42$ \\
        $3$ & $7$ & $2$ & $7$ & $10$ & $12$ & $14$ & $17$ \\
        $3$ & $11$ & $4$ & $9$ & $18$ & $20$ & $24$ & $26$ \\
        \bottomrule
    \end{tabular}
\end{table}
\FloatBarrier

%%%%%%%%%%%%
\section*{A Note on this arXiv Version}
Some minor changes to wording and presentation have been made, and some typos have been corrected. We have added a link \cite{CodeLink} to the SageMath script used for computations, which was not present in the published version of this work. The code has been improved compared to that used for publication, and is significantly more efficient now.

\section*{Acknowledgements}
The authors would like to thank the anonymous referee for comments and corrections which improved the quality of the paper.

\end{document}